\documentclass[11pt]{article}

\usepackage{graphicx}

\newcommand{\bea}{\begin{eqnarray}}
\newcommand{\eea}{\end{eqnarray}}
\newcommand{\be}{\begin{eqpation}}
\newcommand{\ee}{\end{equation}}
\newcommand{\tras}{^{{\mbox{\footnotesize\sc t}}}}

\newcommand\bu {{\bf{u}}}

\newcommand\bx {X}
\newcommand\by {{\bf{y}}}

\newcommand\wf {\widehat{f}}

\newcommand\wtf {\widetilde{f}}

\newcommand\bmu {\mbox{\boldmath $\mu$}}

\def\median{\mathop{\rm median}}

\def\real{\hbox{$\displaystyle I\hskip -3pt R$}}

\def\realito{\tiny{\real}}

\newcommand{\convpp}{ \buildrel{a.s.}\over\longrightarrow}
\newcommand{\convprob}{ \buildrel{p}\over\longrightarrow}
\newcommand{\convdist}{ \buildrel{{\cal D}}\over\longrightarrow}

\newcommand{\MSE}{{\mbox{MSE}}}
\newcommand{\MedSE}{{\mbox{MedSE}}}

\def\dst{\displaystyle}
\def\noi{\noindent}

\def\square{\ifmmode\sqr\else{$\sqr$}\fi}
\def\sqr{\vcenter{
         \hrule height.1mm
         \hbox{\vrule width.1mm height2.2mm\kern2.18mm
\vrule width.1mm}
         \hrule height.1mm}}

\usepackage{color}

\begin{document}

\title{k-Nearest neighbor density estimation on Riemannian Manifolds.}
\author{Guillermo Henry$^1$, Andr\'es Mu\~noz$^{2}$ and Daniela Rodriguez$^1$ \\
{\small \sl $^1$Facultad de Ciencias Exactas y Naturales, Universidad de Buenos Aires and
CONICET, Argentina.}\\
{\small$^2$ Universidad Tecnol\'ogica Nacional and CBC, Universidad de Buenos Aires.}}

\date{}
\maketitle

\begin{abstract}
In this paper, we consider a k-nearest neighbor kernel type estimator when the random variables belong in a Riemannian manifolds. We study asymptotic properties such as the consistency and the asymptotic distribution. A simulation study is also consider to evaluate the performance of the proposal. Finally, to illustrate the potential applications of the proposed estimator, we analyzed two real example where two different manifolds are considered.
\end{abstract}


\section{Introduction}

Let $\bx_1,\dots,\bx_n$ be independent and identically distributed random variables  taking values in $\real^d$ and having density function $f$. A class of estimators of $f$ which has been widely studied since the work of Rosenblatt (1956) and Parzen (1962) has the form
$$f_{n}(x)=\frac{1}{nh^d}{\dst\sum_{j=1}^nK\left(\dst\frac{x-\bx_j}{h}\right)}\;,
$$
where $K(u)$ is a bounded density on $\real^d$ and $h$ is a sequence of positive number such that $h\to 0$ and $nh^d\to \infty$ as $n\to \infty$.

If we apply this estimator to data coming from long tailed distributions, with a small enough $h$ to appropriate for the central part of the distribution, a spurious noise appears in the tails. With a large  value of  $h$ for correctly handling the tails, we can not see the details occurring in the main part of the distribution.  To overcome these defects, adaptive kernel estimators were introduced. For instance, a conceptually similar estimator of $f(x)$ was studied by Wagner (1975) who  defined a general neighbor density estimators by
$$\wf_{n}(x)=\frac{1}{nH_n^d(x)}{\dst\sum_{j=1}^n\dst K\left(\dst\frac{x-\bx_j}{H_n(x)}\right)}\;,
$$
where $H_n(x)$ is the distance between $x$ and the $k$-nearest neighbor of $x$  among $\bx_1,\dots ,\bx_n$, and $k=k_n$ is a sequence of non--random integers such that $\lim_{n\to \infty} k_n=\infty.$ Through  these adaptive bandwidth , the estimation in the point $x$ has the guarantee that to be calculated  using at least  k points of the sample.

However, in many applications, the variables $\bx$ take values on different spaces than $\real^d$. Usually these spaces have a more complicated geometry than the Euclidean space and this has to be taken into account in the analysis of the data. For example if we study the distribution of the stars with luminosity in a given range it is naturally to think that the variables belong to a spherical cylinder ($S^2\times\real$) instead of $\real^4$. If we considerer a region of the planet $M$, then the direction and the velocity of the wind in this region  are points in the tangent bundle of $M$, that is a manifold of dimension $4$. Other examples could be found in image analysis, mechanics, geology and other fields. They include distributions on spheres,  Lie groups, among others, see for example Joshi, et.at. (2007), Goh and Vidal (2008). For this reason, it is interesting to  study an estimation procedure of the density function that take into account a more complex structure of the variables.

 Nonparametric kernel methods for estimating densities of spherical data have been studied by  Hall, et .al (1987) and Bai, et. al. (1988). Pelletier (2005) proposed a family of nonparametric estimators for the density function based on kernel weight when the variables are random object valued in a closed Riemannian  manifold.  The Pelletier's estimators is consistent with the kernel density estimators in the Euclidean case considered by Rosenblatt (1956) and Parzen (1962).

As we comment above, the importance of local adaptive bandwidth is well known in nonparametric statistics and this is even more true with data taking values on complexity space. In this paper, we propose a k-nearest neighbor method when the data takes values on a Riemannian manifolds. The proposal combine the ideas of smoothing in Euclidean spaces with the estimators introduced in Pelletier (2005).

This paper is organized as follows. Section \ref{preliminares} contains a brief summary of the  necessaries concepts of Riemannian geometry. In Section \ref{estimador}, we introduce a k-nearest neighbor estimators on Riemannian manifolds. Uniform consistency of the estimators is derived in Section \ref{consistencia}, while  in Section \ref{distribucion} the asymptotic distribution is obtained under regular assumptions.  Section \ref{simulaciones} contains a Monte Carlo  study designed to evaluate the proposed estimators. Finally, Section \ref{ejemplos} presents two example using real data. Proofs are given in the Appendix.

\section{Preliminaries and the estimator}{\label{preliminares}}
Let $({M},g)$ be a $d-$dimensional   Riemannian manifold without boundary. We denote by $d_g$ the distance induced by the metric  $g$. With $B_s(p)$ we denote a normal ball with radius $s$ centered at $p$. The injectivity radius of $(M,g)$ is given  by
 $inj_g M=\dst\inf_{p\in M}\sup \{s \in \real >0: B_s(p) \mbox{ is a normal ball} \}$.  Is easy to see that  a compact Riemannian manifold  has strictly positive injectivity radius. For example, it is not difficult to see that the $d$-dimensional sphere $S^d$ endowed with the metric induced by the canonical metric $g_0$ of $R^{d+1}$ has injectivity radius equal to $\pi$. If $N$ is a proper submanifold of the same dimension than $(M,g)$, then $inj_{g|_{N}} N=0$. The Euclidean space or the hyperbolic space have infinity injectivity radius. Moreover, a complete and simply connected Riemannian manifold with non positive sectional curvature has also this property.

 Throughout this paper, we will assume that $(M,g)$ is a complete Riemannian manifold, i.e. $(M,d_g)$ is a complete  metric space. Also we will consider that $inj_g M$ is strictly positive.  This restriction will be clear in the Section \ref{estimador} when we define the estimator. For standard result on differential and Riemannian geometry we refer to the reader to     Boothby (1975), Besse (1978), Do Carmo (1988) and Gallot, Hulin and Lafontaine (2004).

Let $p\in M$, we denote with $0_p$ and $T_pM$ the null tangent vector and the tangent space of $M$ at $p$. Let $B_s(p)$ be a normal ball centered at $p$. Then $B_s(0_p)=exp_p^{-1}(B_s(p))$ is
an open neighborhood of $0_p$ in $T_pM$ and so it has a natural
structure of differential manifold. We are going to consider the
Riemannian metrics $g^{\,\prime}$ and $g^{\,\prime\,\prime}$ in
$B_s(0_p)$, where $g^{\,\prime}=exp_p^*(g)$ is the pullback of $g$
by the exponential map and $g^{\,\prime\,\prime}$ is the canonical
metric induced by $g_p$ in $B_s(0_p)$. Let $w\in B_s(0_p)$, and
 $(\bar{U},\bar{\psi})$ be a chart of $B_s(0_p)$ such that $w\in
\bar{U}$. We note by $\{{\partial}/{\partial
\bar{\psi}_1}|_w,\dots,{\partial}/{\partial
\bar{\psi}_d}|_w\}$ the  tangent vectors induced by $(\bar U,\psi)$. Consider the matricial function with entries $(i, j)$ are given by $g^{\,\prime}\left(\left({\partial}/{\partial
\bar{\psi}_i}|_w\right),\left({\partial}/{\partial
\bar{\psi}_j}|_w\right)\right)$. The volumes of the parallelepiped spanned by
$\{\left({\partial}/{\partial
\bar{\psi}_1}|_w\right),\dots,\left({\partial}/{\partial
\bar{\psi}_d}|_w\right)\}$ with respect to the metrics
$g^{\,\prime}$ and $g^{\,\prime\,\prime}$ are given by ${|\det
g^{\,\prime}\left(\left({\partial}/{\partial
\bar{\psi}_i}|_w\right),\left({\partial}/{\partial
\bar{\psi}_j}|_w\right)\right)|^{1/2}}$ and ${|\det
g^{\,\prime\,\prime}\left(\left({\partial}/{\partial
\bar{\psi}_i}|_w\right),\left({\partial}/{\partial
\bar{\psi}_j}|_w\right)\right)|^{1/2}}$ respectively. The quotient
between these two volumes is independent of the selected chart.
So, given $q\in B_s(p)$,  if $w=exp_p^{-1}(q)\in B_s(0_p)$ we can
define the volume density
function, $\theta_p(q)$, on $(M,g)$  as
$$\theta_p(q)=
\frac{{|\det g^{ \,\prime}\left(\left({\partial}/{\partial
\bar{\psi}_i}|_w\right) ,\left({\partial}/{\partial
\bar{\psi}_j}|_w\right)\right)|}^{1/2}} {{|\det
g^{\,\prime\,\prime}\left(\left({\partial}/{\partial
\bar{\psi}_i}|_w\right),\left({\partial}/{\partial
\bar{\psi}_j}|_w\right)\right)|}^{ 1/2} }$$
 for any chart
$(\bar{U},\bar{\psi})$ of $B_s(0_p)$ that contains
$w=exp^{-1}_p(q)$. For instance, if we consider a normal coordinate system  $(U,\psi)$  induced by an orthonormal basis
$\{v_1,\dots,v_d\}$ of $T_pM$ then $\theta_p(q)$ is the function of the volume element $d\nu_g$   in the local expression with respect to chart $(U,\psi)$ evaluated at $q$, i.e.
$$\theta_p(q)={\left|\det g_q\left(\frac{\partial}{\partial
\psi_i}\Big|_q,\frac{\partial}{\partial \psi_j}\Big|_q
\right)\right|}^{\frac 12}\ ,$$ where
$\frac{\partial}{\partial\psi_i}|_q=D_{\alpha_i(0)}exp_p(\dot{\alpha}_i(0))$ with $\alpha_i(t)=exp_p^{-1}(q)+tv_i$ for $q\in U$.
Note that the volume density function $\theta_p(q)$ is not defined for all the pairs $p$ and $q$ in $M$, but it is if $d_g(p,q)<inj_g M$.

We finish the section showing some examples of the density function:

\begin{enumerate}
\item[i)] In the case of  $(\real^d,g_0)$ the density function is  $\theta_{p}(q)=1$ for all $(p,q)\in\real^d\times\real^d$.
\item[ii)] In the  $2$-dimensional sphere of radius $R$,  the volume density is $$\theta_{p_1}(p_2)=R\;\frac{|sen(d_g(p_1,p_2)/R)|}{d_g(p_1,p_2)}\quad \mbox{ if } p_2\neq p_1, -p_1\  \mbox{ and } \theta_{p_1}(p_1)=1.$$
    where $d_g$ induced is given by $$d_g(p_1,p_2)=R \arccos(\frac{{< \; p_1,\;p_2>}}{R^2}).$$

\item[iii)] In the case of the cylinder of radius 1 ${\cal{C}}_1$ endowed with the metric induced by the canonical metric of $\real^3$,  $\theta_{p_1}(p_2)=1$ for all $(p_1,p_2)\in $${\cal{C}}_1\times{\cal{C}}_1 $, and the distance induced is given by $d_g(p_1,p_2)=d_2((r_1,s_1),(r_2,s_2))$ if $d_2((r_1,s_1),(r_2,s_2))<\pi$, where $d_2$ is the Euclidean distance of $\real^2$ and $p_i=(\cos(r_i),\sin(r_i),s_i)$ for $i=1,2$.

\end{enumerate}
See also, Besse (1978) and Pennec (2006) for a discussion on the volume density function.

\subsection{The estimator}\label{estimador}

Consider a probability distribution with a density $f$ on a $d-$dimensional Riemannian manifold $(M,g)$. Let $\bx_1,\cdots, \bx_n$ be i.i.d  random object takes values on $M$  with density $f$.
A natural extension of the estimator proposed by Wagner (1975) in the context of a Riemannian manifold is to consider the following  estimator
$$\wf_{n}(p)=\frac{1}{nH_n^d(p)}{\dst\sum_{j=1}^n\dst\frac{1}{\theta_{\bx_j}(p)}K\left(\dst\frac{d_g(p,\bx_j)}{H_n(p)}\right)}\;,
$$
where $K:\real \to \real$ is a non-negative function with compact support, $\theta_p(q)$ denotes the volume density
function on $(M,g)$ and   $H_n(p)$ is the distance $d_g$ between $p$ and the $k$-nearest neighbor of $p$  among $\bx_1,\dots ,\bx_n$, and $k=k_n$ is a sequence of non--random integers such that $\lim_{n\to \infty} k_n=\infty.$

As we mention above, the volume density function is not defined for all $p$ and $q$. Therefore, in order to guarantee the well definition of the estimator we consider a modification of the proposed estimator. Using the fact that the kernel $K$ has compact support,    we consider as  bandwidth  $\zeta_n(p)=\min\{H_n(p),inj_g M\}$ instead of   $H_n(p)$. Thus, the kernel only consi-ders the points $\bx_i$ such that $d_g(\bx_i,p)\leq \zeta_n(p)$ that are smaller than $inj_g M$ and for these points, the volume density function is well defined. Hence, the k-nearest neighbor kernel type estimator is defined as follows,
\begin{eqnarray}\label{estim}
\wf_{n}(p)=\frac{1}{n\zeta_n^d(p)}{\dst\sum_{j=1}^n\dst\frac{1}{\theta_{\bx_j}(p)}K\left(\dst\frac{d_g(p,\bx_j)}{\zeta_n(p)}\right)}\;,
\end{eqnarray}
where $\zeta_n(p)=\min\{H_n(p),inj_g M\}$.

\noi \textbf{Remark \ref{estimador}.1.}  If $(M,g)$ is a compact Riemannian manifold and its sectional curvature is not bigger than $a>0$, then we know by the Lemma of Klingerberg  (see Gallot, Hulin, Lafontaine (2004)) that $inj_g M\geq min\{\pi/\sqrt{a},l/2\}$ where $l$ is the length of the shortest closed geodesic in $(M,g)$.

\section{Asymptotic results}
Denote by $C^{\ell}(U)$ the set of  $\ell$ times continuously differentiable functions from $U$ to $\real$ where $U$ is an open set of $M$.
We assume that the measure induced by the probability $P$ and by $\bx$ is absolutely continuous with respect to the Riemannian volume measure $d\nu_g$, and we denote by $f$  its density  on $M$  with respect to  $d\nu_g$. More precisely, let ${\cal B}(M)$ be the Borel $\sigma-$field of $M$ (the $\sigma-$field generated by the class of open sets of $M$). The random variable $\bx$ has a probability density function $f$, i.e. if $\chi\in{\cal B}(M)$, $P(\bx^{-1}(\chi))=\int_{\chi} f d\nu_g.$

\subsection{Uniform Consistency}\label{consistencia}
We will consider the following set of assumptions in order  to derive the strong consistency results of the estimate $\wf_n(p)$ defined in (\ref{estim}).

\begin{enumerate}
\item[$H1.$] Let ${M}_0$ be a compact set on ${M}$ such that:
\begin{enumerate}
\item[i)] $f$ is a bounded function such that $\inf_{p\in M_0}f(p)=A>0.$
\item[ii)] $\inf_{p,q\in M_0}\theta_p(q)=B>0.$
\end{enumerate}
\item[$H2.$] For any open set $U_0$ of $M_0$ such that $M_0\subset U_0$,
$f$ is of class $C^2$ on $U_0$.
\item[$H3.$] The sequence $k_n$ is such that $k_n\to \infty$, $\frac{k_n}{n}\to 0$ and $\frac{k_n}{\log n}\to \infty$  as $n\to \infty$.
\item[$H4.$] $K: \real \to \real$ is a bounded
nonnegative Lipschitz function of order one, with compact support $[0,1]$   satisfying: $\int_{\realito^d}
K(\|\bu\|)d\bu=1$, $\int_{\realito^d} \bu
K(\|\bu\|)d\bu=\bf{0}$ and \linebreak $0<\int_{\realito^d}
\|\bu\|^2K(\|\bu\|)d\bu<\infty$.
\item[$H5.$] The kernel $K(u)$ verifies $K(uz)\geq K(z)$ for all $u\in (0,1)$.
\end{enumerate}

\noi \textbf{Remark \ref{consistencia}.1.} The fact that $\theta_p(p)=1$ for all $p\in M$ guarantees that $H1$ ii)  holds. The assumption $H3$ is usual when dealing with nearest neighbor  and the assumption $H4$   is  standard  when dealing with kernel estimators.

\noi \textbf{Theorem \ref{consistencia}.2.} Assume that  $H1$ to  $H5$ holds, then we have that
$$
 \dst\sup_{p\in {{ M}_0}}|\wf_n(p)-f(p)|\convpp 0.
$$

\subsection{Asymptotic normality}\label{distribucion}
To derive the asymptotic distribution of the regression parameter estimates we will need two additional assumptions. We will denote with ${\cal V}_r$ the Euclidean ball of radius $r$ centered at the origin and with  $\lambda({\cal V}_r)$  its  Lebesgue measure.

\begin{enumerate}
\item[$H5.$] $f(p)>0$, $f\in C^{2}(V)$ with $V\subset M$ an open neighborhood of $M$ and the second derivative of $f$ is bounded.
\item[$H6.$] The sequence $k_n$ is such that $k_n\to \infty$, $k_n/n\to 0$ as $n\to \infty$ and  there exists $0\leq \beta <\infty$ such that $\sqrt{k_n n^{-4/(d+4)}}\to \beta$ as $n\to \infty.$
\item[$H7.$] The kernel verifies
 \begin{enumerate}
 \item[i)] $\int K_1(\|\bu\|)\|\bu\|^2d\bu<\infty$  as $s\to\infty$ where $K_1(\bu)=K'(\|\bu\|)\|\bu\|$.
 \item[ii)] $\|\bu\|^{d+1}K_2(\bu)\to 0$ as $\|\bu\|\to\infty$   where $K_2(\bu)=K''(\|\bu\|)\|\bu\|^2-K_1(\bu)$
\end{enumerate}
\end{enumerate}

\noi \textbf{Remark \ref{distribucion}.1.}  Note that $div(K(\|\bu\|)\bu)=K'(\|\bu\|)\|\bu\|+d\;K(\|\bu\|)$, then using the divergence Theorem, we get that $\int K'(\|\bu\|)\|\bu\|d\bu =\int_{\|\bu\|=1} K(\|\bu\|)\bu \frac{\bu }{\|\bu\|} d\bu - d\; \int K(\|\bu\|) d\bu$. Thus, the fact that $K$ has compact support in $[0,1]$  implies that $\int K_1(\bu) d\bu=-d.$

On the other hand, note that $\nabla(K(\|\bu\|)\|\bu\|^2)=K_1(\|\bu\|)\bu+2K(\|\bu\|)\bu$ and by $H4$
we get that $\int K_1(\|\bu\|)\bu d\bu={\bf 0}$.

\noi \textbf{Theorem \ref{distribucion}.2.} Assume $H4$ to $H7$. Then we have that
$$
\sqrt{k_n}(\wf_n(p)-f(p))\convdist {\cal N}(b(p),\sigma^2(p))
$$
with
$$
b(p)= \frac 12 \frac{\beta^{\frac{d+4}{d}}}{{(f(p)\lambda({\cal V}_1))}^{\frac 2d}}  \dst\int_{{\cal V}_1}K(\|\bu\|)u_1^2 \;d\bu\; \dst\sum_{i=1}^d\frac{\partial f\circ \psi^{-1}}{\partial u_i\partial u_i}|_{u=0}\;
$$
and
$$
\sigma^2(p)=\lambda({\cal V}_1)f^2(p)\int_{{\cal V}_1} K^2(\|\bu\|)d\bu
$$
where $\bu=(u_1,\dots,u_d)$  and $(B_{h}(p),\psi)$ is any  normal coordinate system.

In order to derive the asymptotic distribution of $\wf_n(p)$, we will study  the asymptotic behavior of ${h^d_n}/{\zeta^d_n(p)}$ where  $h_n^d={k_n}/{(n f(p)\lambda({\cal V}_1))}$. Note that if we consider $\wtf_n(p)={k_n}/{(n \zeta^d_n(p) \lambda({\cal V}_1))}$,  $\wtf_n(p)$ is a consistent estimator of $f(p)$ (see  the proof of Theorem \ref{consistencia}.2.). The next Theorem states  that this estimator is also asymptotically  normally distributed as in the Euclidean case.
\vspace{0.5cm}

\noi \textbf{Theorem \ref{distribucion}.3.}  \textsl{Assume $H4$ to $H6$, and
 let $h_n^d={k_n}/({n f(p)\lambda({\cal V}_1)})$. Then we have that
$$
\sqrt{k_n}\left(\frac{h^d_n}{\zeta^d_n(p)}-1\right)\convdist N(b_1(p),1)
$$
with
$$b_1(p)=\left(\frac{\beta^{\frac{d+4}2}}{f(p)\mu({\cal V}_1)}\right)^{\frac 2d} \left\{\frac{\tau}{6d+12}+ \frac{\int_{{\cal V}_1}u_1^2\,d\bu\;L_1(p)}{f(p)\mu({\cal V}_1)}\right\}$$ where $\bu=(u_1, \dots, u_d)$, $\tau $ is the scalar curvature of $(M,g)$, i.e. the trace of the Ricci tensor,
$$L_1(p)=\dst\sum_{i=1}^d\left(\displaystyle\frac{\partial^2f \circ
\psi^{-1}}{\partial u_iu_i }\Big|_{u=0}+\displaystyle\frac{\partial f \circ
\psi^{-1}}{\partial u_i }\Big|_{u=0} \displaystyle\frac{\partial \theta_p \circ
\psi^{-1}}{\partial u_i }\Big|_{u=0} \right)$$
  and $(B_h(p), \psi)$ is any normal coordinate system.}

\section{Simulations}\label{simulaciones}

This section contains the results of a simulation study designed to evaluate the performance of the estimator defined in the Section \ref{estimador}.
We consider three models in two different Riemannian manifolds, the sphere and the cylinder endowed with the metric induced by the canonical metric of $\real^3$. We performed 1000 replications of independent samples of size $n=200$ according to the following models:

\begin{enumerate}
\item[] \bf Model 1 (in the sphere)\rm: The random variables $\bx_i$ for $1\leq i\leq n$ are i.i.d.  Von Mises distribution  $ VM(\mu,\kappa)$ i.e.  $$f_{\mu,\kappa}(\bx)=\left(\frac{k}{2}\right)^{1/2} I_{1/2}(\kappa) \exp\{\kappa\,\bx\tras \bmu \},$$ with $\bmu$  is the mean parameter, $\kappa>0$ is the concentration parameter and $I_{1/2}(\kappa)=\left(\frac{\kappa\pi}{2}\right) \sinh(\kappa)$ stands for the modified Bessel function. This model has many important applications, as described in \cite{js} and \cite{mj}. We generate a random sample follows a Von Mises distribution with mean $(0,0,1)$ and concentration parameter $3$.
\item[] \bf Model 2 (in the sphere)\rm: We simulate i.i.d. random variables $Z_i$ for $1\leq i\leq n$ following a multivariate normal distribution of dimension $3$, with mean $(0,0,0)$ and covariance matrix  equals to the identity. We define $\bx_i=\frac{Z_i}{\|Z_i\|}$ for $1\leq i\leq n$,  therefore the variables $\bx_i$ follow an uniform distribution in the two dimensional sphere.
\item[]\bf Model 3 (in the cylinder)\rm:
 We consider  random variables $\bx_i=(\by_i,t_i)$ taking values in the cylinder $S^1\times\real$. We generated the model proposed by Mardia and Sutton (1978) where,
$$
\by_i=(\cos(\theta_i),\sin(\theta_i)) \sim VM((-1,0),5)$$ $$ t_i|\by_i\sim N(1+2\sqrt{5}\cos(\theta_i),1).
$$
Some examples of variables with this distribution can be found in Mardia and Sutton (1978).
\end{enumerate}

In all cases, for smoothing procedure, the kernel was taken as the quadratic kernel $K(t)=( {15}/{16}) (1-t^2)^2 I(|x|<1)$. We have considered a grid of equidistant values of $k$ between $5$ and $150$ of length $20$.

 To study the performance of the estimators of the density function $f$, denoted by  $\wf_n$, we have considered the mean square error  ($\MSE$) and the median square  error ($\MedSE$), i.e,
$$
\MSE(\wf_n)=\frac 1n \sum_{i=1}^n[\wf_n(\bx_i)-f(\bx_i)]^2\;.
$$

$$
\MedSE(\wf_n)=\median|\wf_n(\bx_i)-f(\bx_i)|^2\;.
$$
Figure 1  gives the values of the  $\MSE$ and ${\MedSE}$ of $\wf_n$ in  the sphere model considering  different  numbers of neighbors, while Figure 2 shows the cylinder model.
The simulation study confirms the good behavior of $k-$nearest neighbor estimators,  under the different models considered. In all cases, the estimators are stable under large numbers of neighbors. However, as expected,  the estimators using a small number of neighbors  have a poor  behavior, because in the neighborhood of each point there is a small number of samples.

\begin{center}
\vspace{-2cm}
\includegraphics[scale=0.34,angle=90]{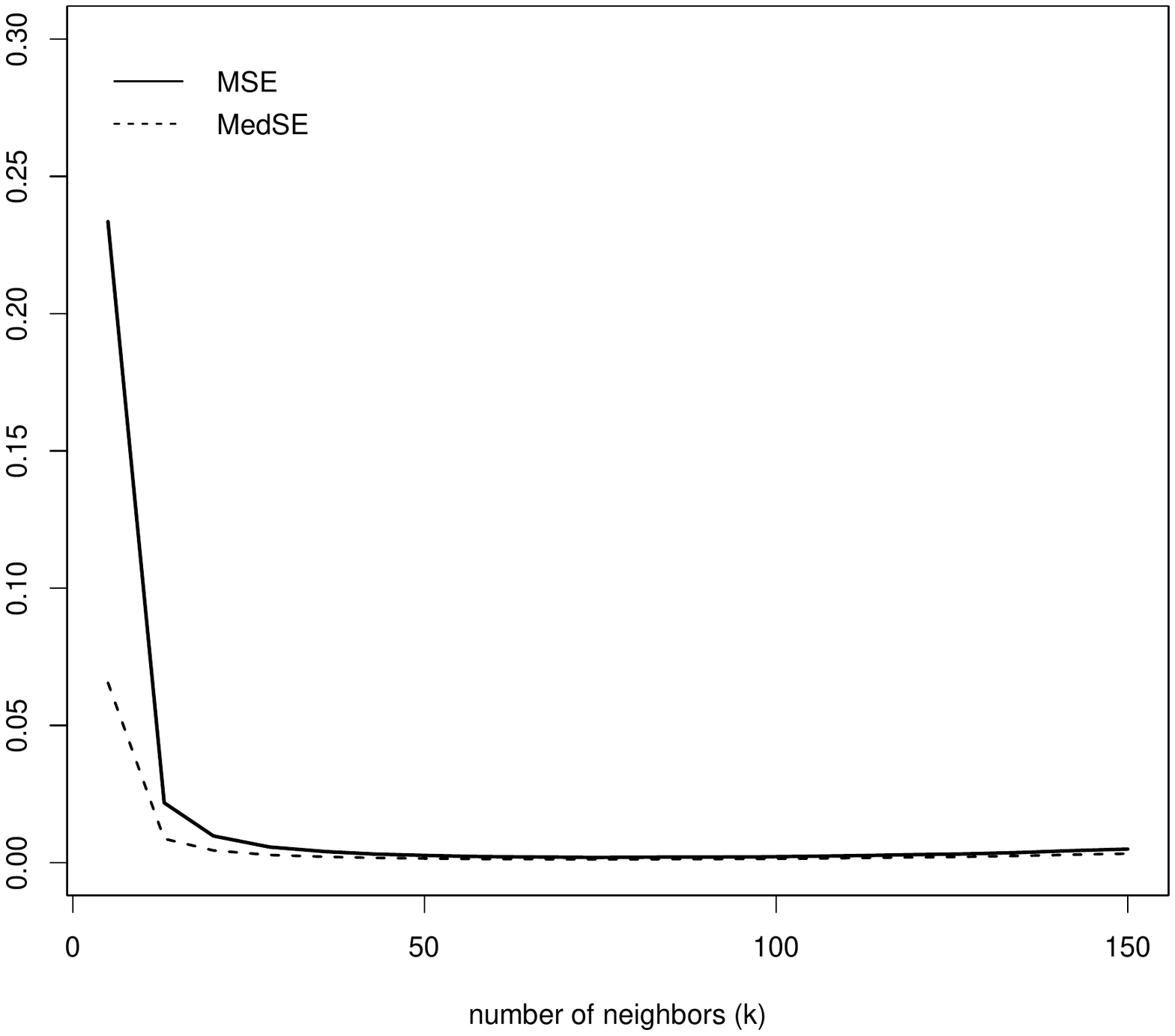}\hskip0in \includegraphics[scale=0.34,angle=90]{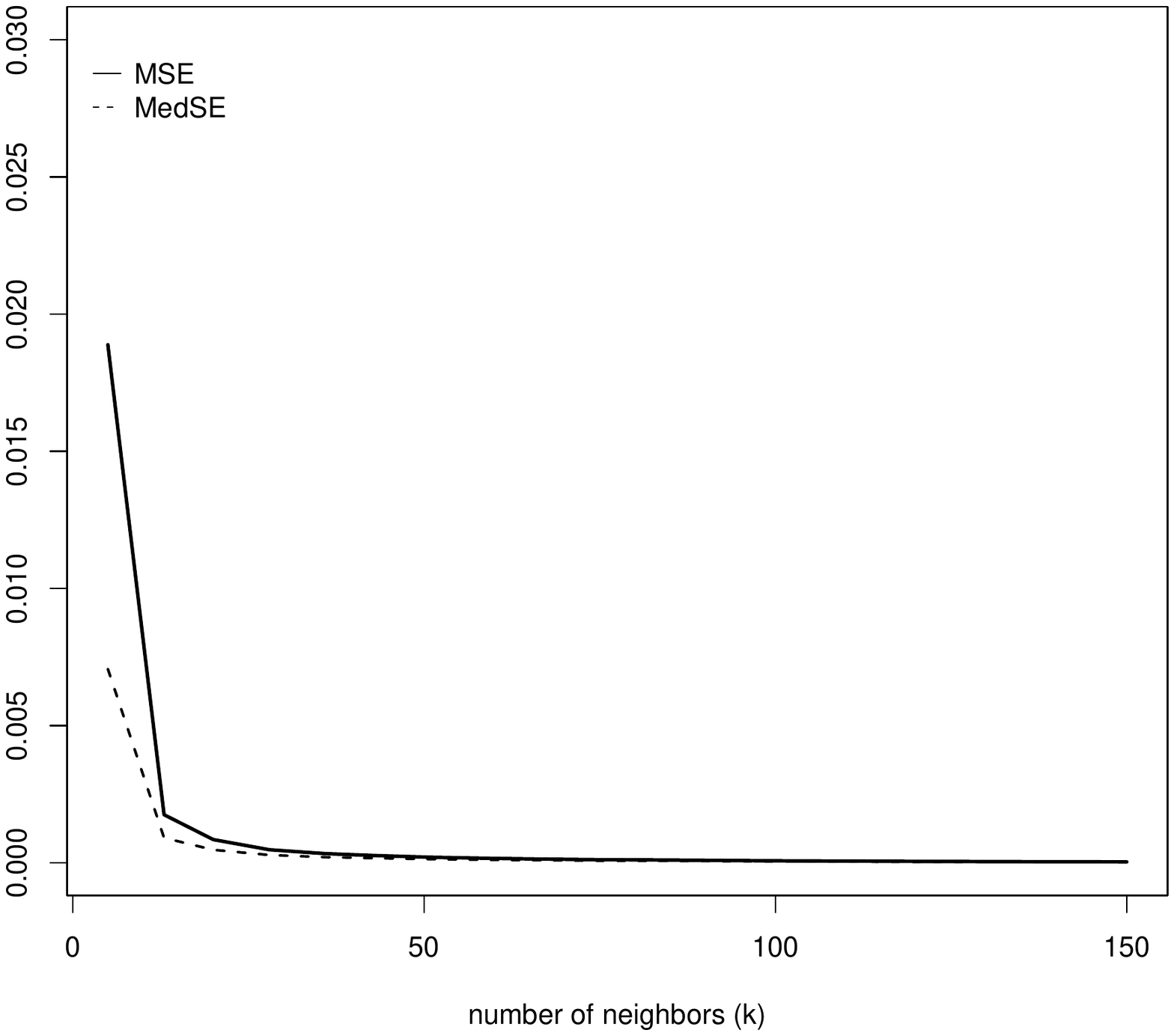}\\
\vspace{-2cm}
$\quad\quad\quad\quad\quad a)$\hskip2.5in $b)\quad\quad\quad\quad\quad\quad\quad$\\
\end{center}
\vspace{-0.5cm}
\small Figure 1: The nonparametric density estimator using different numbers of neighbor, a) the Von Mises model and  b) the uniform model.\normalsize

\begin{center}
\vspace{-2.5cm}
\includegraphics[scale=0.34]{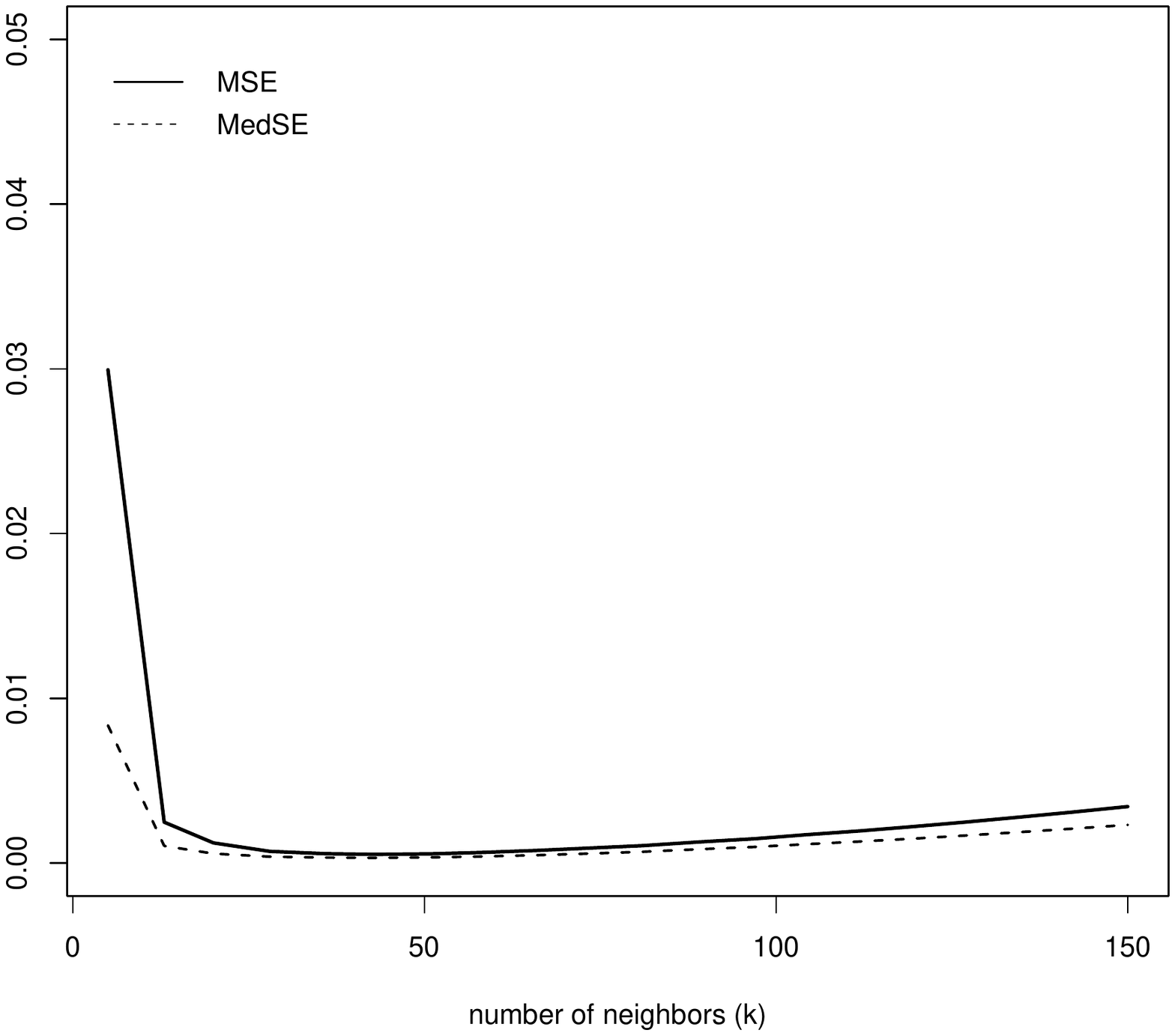}\\
\vspace{-2cm}
$c)$\\
\end{center}
\small Figure 2: The nonparametric density estimator using different numbers of neighbor in the cylinder.\normalsize

\section{Real Example}\label{ejemplos}

\subsection{Paleomagnetic data}\label{volvanes}

It is well know the need for statistical analysis of paleomagnetic data. Since the work developed by Fisher (1953),  the study  of parametric families was  considered a principal tools to analyze  and quantify this type of data (see  Cox and Doell (1960), Butler (1992) and Love and Constable (2003)). In particular, our proposal allows to explore the nature of directional dataset that include paleomagnetic data without make any parametric assumptions.

In order to illustrate the k-nearest neighbor kernel type estimator  on the two-dimensional sphere, we  illustrate the estimator using  a paleomagnetic data set studied by Fisher, Lewis, and  Embleton (1987).  The data set consist in  $n=107$ sites from specimens of Precambrian volcanos whit measurements of magnetic remanence. The data set contains two variables corresponding to the directional component on a longitude scale, and the directional component on a latitude scale. The original data set is available in the library \sf sm \rm of \sf R \rm statistical package.

To calculate the estimators the volume density function and the geodesic distance  were taken as in the Section \ref{preliminares} and we considered the quadratic kernel $K(t)=( {15}/{16}) (1-t^2)^2 I(|x|<1)$. In order to analyzed the sensitivity of the results with respect to the number  of neighbors, we plot  the estimator using  different bandwidths. The results are shown in the Figure 3.

The real data  was  plotted in blue and with a large radius in order to obtain a better visualization. The  Equator line, the  Greenwich  meridian and  a second meridian are in gray while the north and south pols are denoted with the capital letter  N and S respectively. The levels of concentration  of measurements of magnetic remanence  can be found in yellow for high levels  and  in red for lowest density levels.  Also, the levels   of concentration  of measurements of magnetic remanence  was illustrated with relief on  the sphere that allow to emphasize  high density levels and the form of the density function.

As in the Euclidean case large number of neighbors produce estimators with small variance but high bias, while small
values produce more wiggly estimators.  This fact shows the need of the implementation of a  method to select the adequate bandwidth for this estimators.  However, this require further careful investigation and are beyond the scope of this paper.

\begin{center}
\includegraphics[scale=0.8]{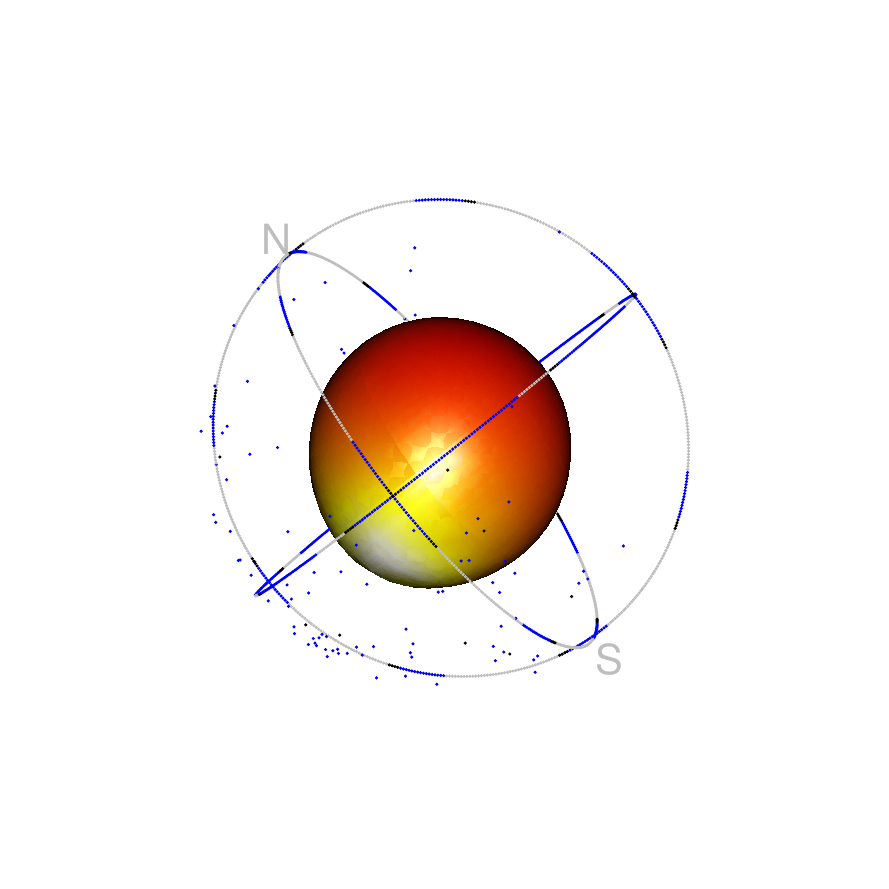}\hspace{-0.5cm} \includegraphics[scale=0.8]{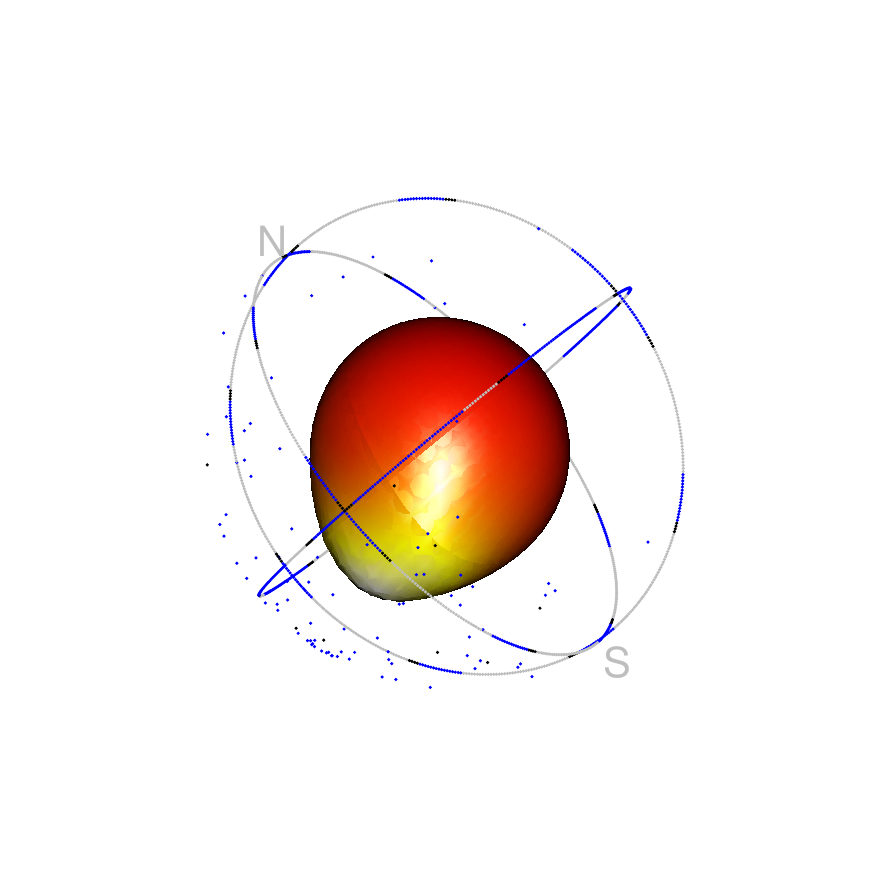}\\
\vspace{-1.5cm}
$a)$\hskip2.8in $b)$\\
\vspace{-0.8cm}
{\includegraphics[scale=0.8]{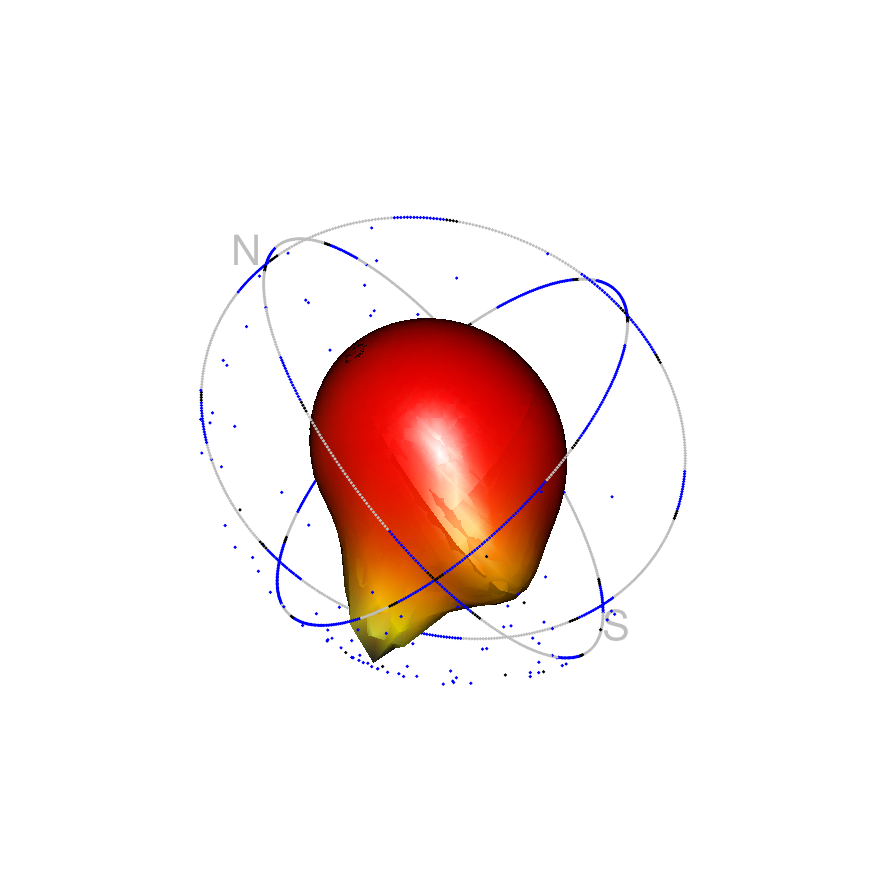}}\hspace{-0.5cm} {\includegraphics[scale=0.75]{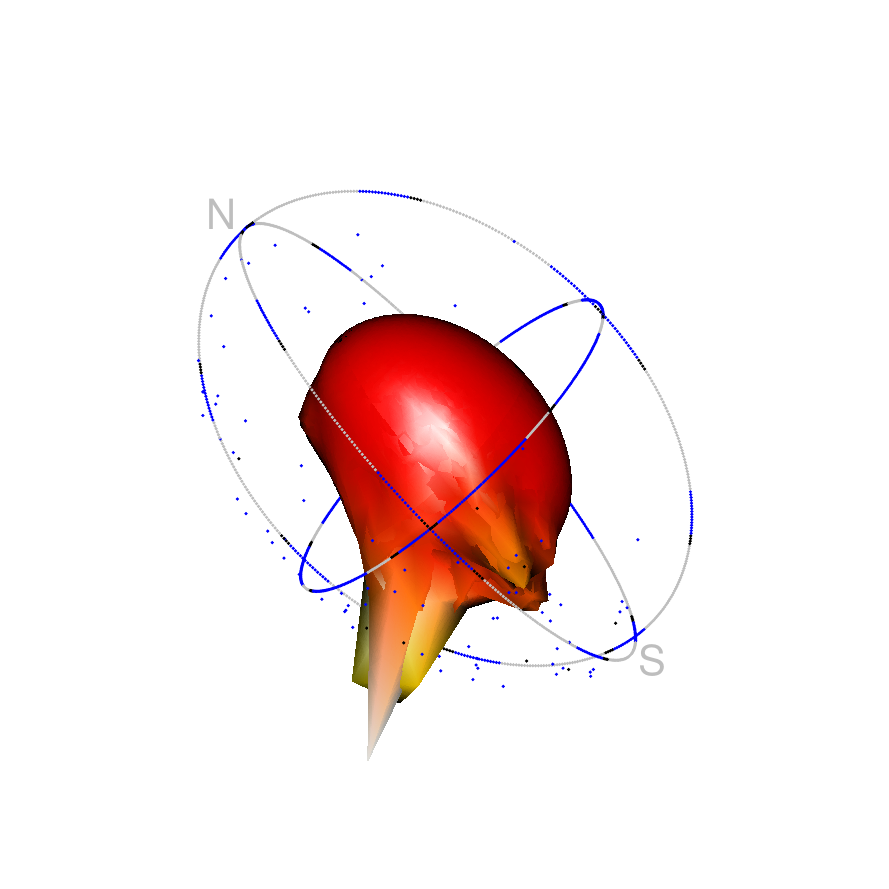}}\\
\vspace{-1.5cm}
$c)$\hskip2.8in $d)$\\
\end{center}
\small Figure 3: The nonparametric density estimator using different number of neighbors, a) k=75, b)k=50,  c) k=25 and d) k=10.

\subsection{Meteorological data}

In this Section, we consider a real data set collected in the  meteorological station \lq\lq Ag\"uita de Perdiz\rq\rq $\;$  that is located in Viedma, province of R\'io Negro, Argentine. The dataset consists in the wind direction and temperature during January 2011 and  contains $1326$ measurements that were registered  with a frequency of   thirty minutes. We note that the considered variables  belong to a cylinder with radius 1.

As in the previous Section, we consider the quadratic kernel and  we took the density function  and the geodesic distance as in Section \ref{preliminares}. Figure 4 shows the result of the estimators, the color and form of the graphic was constructed as in the previous example.

It is important to remark that  the measurement devices of  wind direction  not present a sufficient precision to avoid  repeated data. Therefore, we consider the  proposal given  in Garc\'ia-Portugu\'es, et.al. (2011) to solve this problem. The proposal consists in  perturb the repeated data as follows
$\widetilde{r}_i = r_i + \xi \varepsilon_i,$
where $r_i$ denote the wind direction measurements and $\varepsilon_i,$  for $i = 1,\dots, n$  were independently
generated from a von Mises distribution with $\mu =(1,0)$ and $\kappa = 1$. The selection of the
perturbation scale $\xi$ was taken  $\xi = n^{-1/5}$  as in Garc\'ia-Portugu\'es, et.al. (2011) where in this case $n=1326$

The work of Garc\'ia-Portugu\'es, et.al. (2011)  contains other nice real example where   the proposed  estimator can be apply. They considered a naive density estimator applied to  wind directions and  SO2 concentrations, that allow you explore high  levels of contamination.

\begin{center}
\vspace{-1cm}
\includegraphics[scale=0.8]{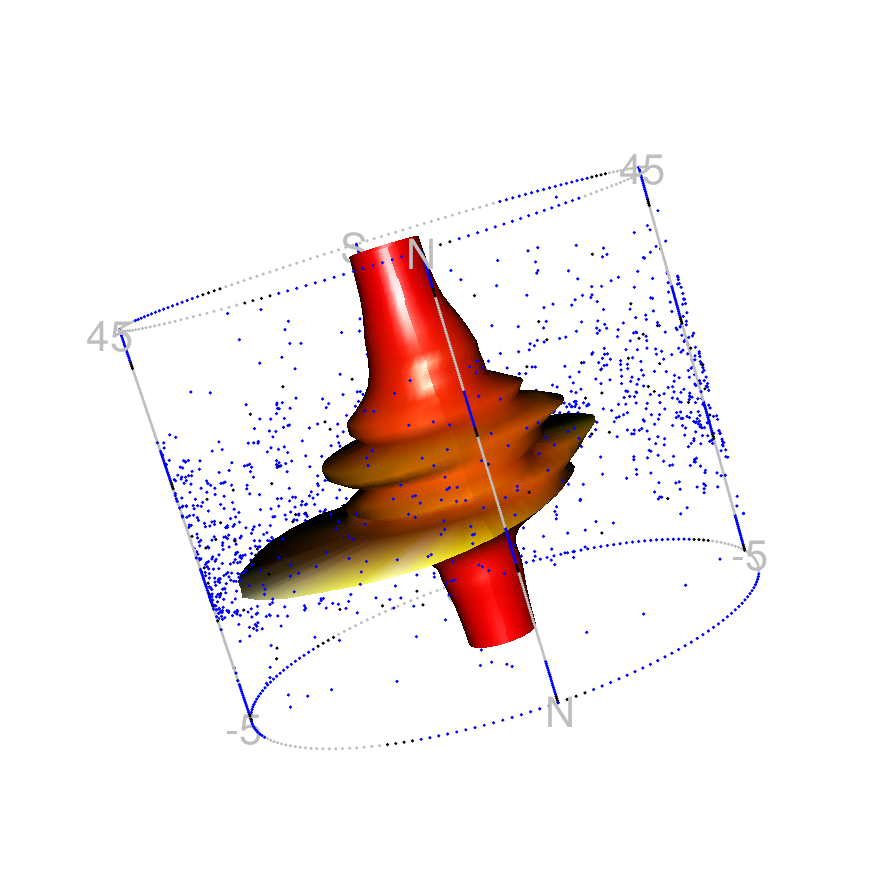}\hspace{0.5cm} \includegraphics[scale=0.8]{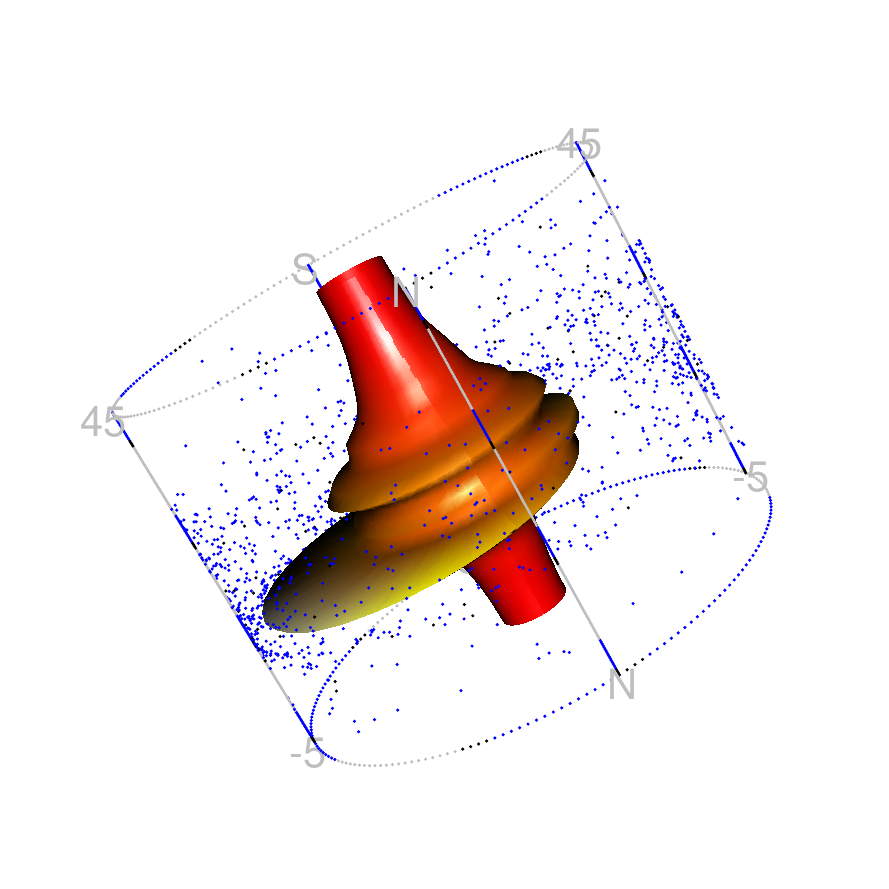}\\
\vspace{-0.8cm}
$a)$\hskip3in $b)$\\
\vspace{-1cm}
{\includegraphics[scale=0.8]{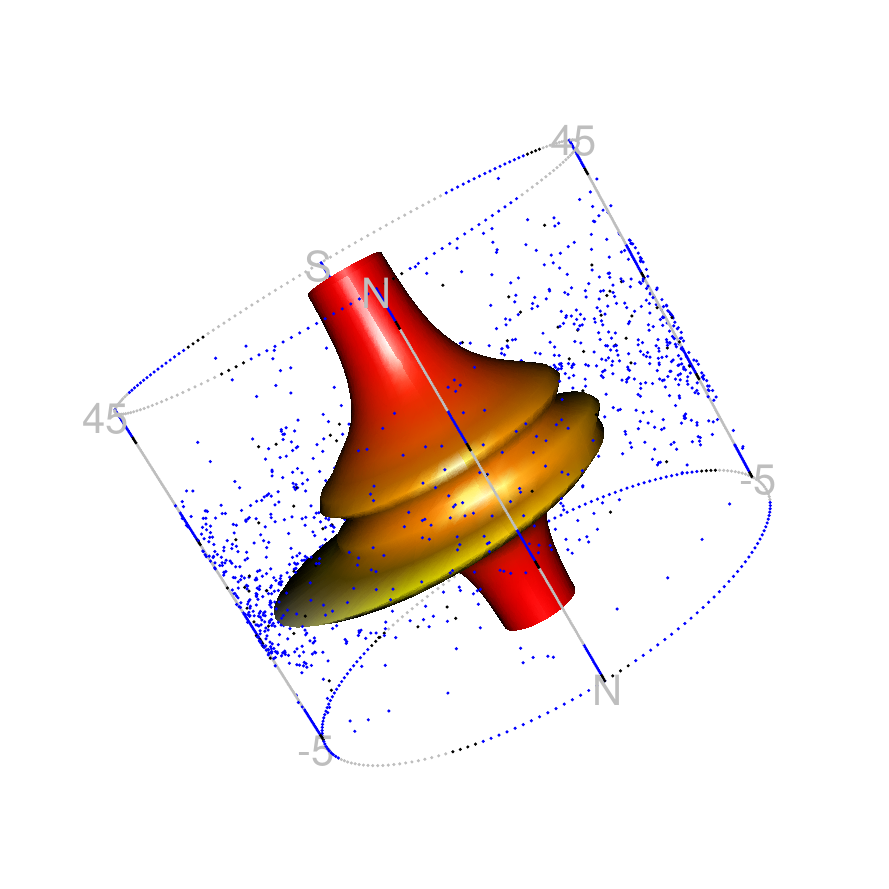}}\hspace{0.5cm} {\includegraphics[scale=0.8]{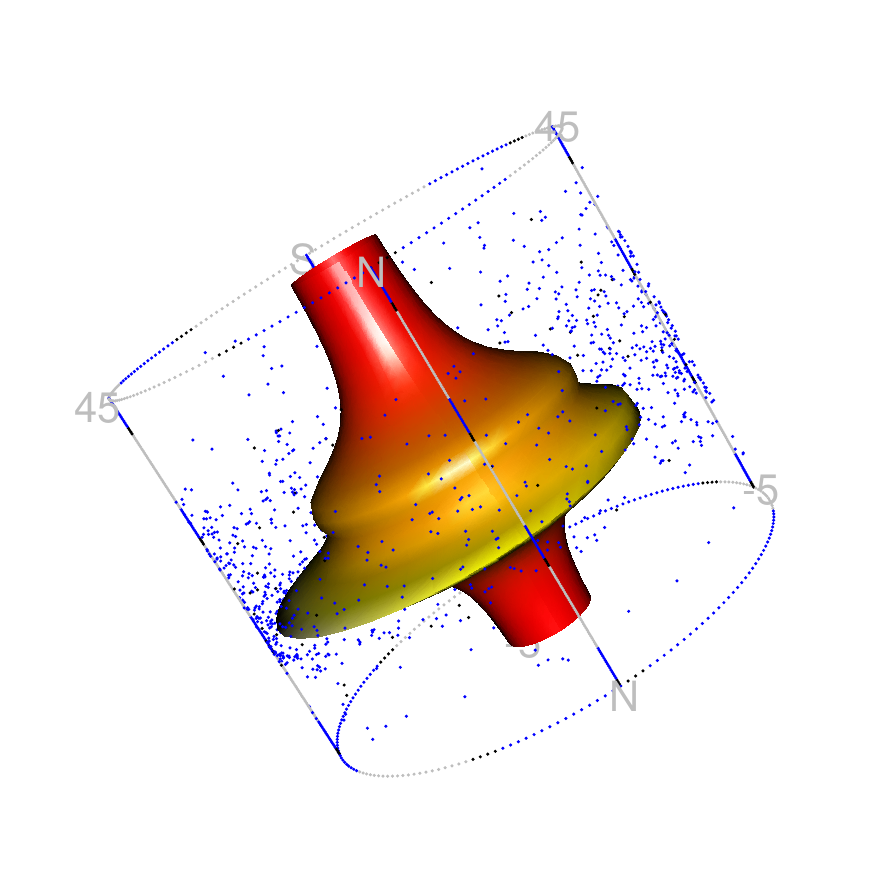}}\\
\vspace{-0.8cm}
$c)$\hskip2.8in $d)$\\
\end{center}
\small Figure 4: The nonparametric density estimator using different number of neighbors, a) k=75, b)k=150,  c) k=300 and d) k=400.

In Figure 4, we can see that the lowest temperature are more probable when the wind  comes from the East direction. However, the highest temperature does not seem to have  correlation with the wind direction. Also, note that in Figure 4, we can see two mode  corresponding  to the minimum and maximum daily of the temperature.

These examples show the usefulness  of the proposed estimator for the analysis and exploration of these type of dataset.


\thispagestyle{empty} 
\appendix
\section*{Appendix}\label{proofs}

\subsection*{Proof of Theorem \ref{consistencia}.2. }
Let $$
f_n(p,\delta_n)=\frac{1}{n\delta_n^d} \sum_{i=1}^n \frac{1}{\theta_{\bx_i}(p)}K\left(\frac{d_g(p,\bx_i)}{\delta_n}\right)\;.
$$
Note that if $\delta_n=\delta_n(p)$ verifies $\delta_{1n}\leq\delta_n(p)\leq\delta_{2n}$ for all $p\in M_0$ where  $\delta_{1n}$ and $\delta_{2n}$ satisfy $\delta_{in}\to 0$ and $\dst\frac{n\delta_{in}^d}{\log n}\to \infty $ as $n\to \infty$ for $i=1,2$ then by Theorem 3.2 in Henry and Rodriguez (2009) we have that
\begin{eqnarray}
 \dst\sup_{p\in {{M}_0}}|f_n(p,\delta_n)-f(p)|\convpp 0\label{sinvecinos}
\end{eqnarray}

For each $0<\beta<1$ we define,
$$
f^-_n(p,\beta)=\frac{1}{nD^+_n(\beta)^d} \sum_{i=1}^n \frac{1}{\theta_{\bx_i}(p)}K\left(\frac{d_g(p,\bx_i)}{D^-_n(\beta)}\right)=f^-_n(p,D^-_n(\beta)^d)\frac{D^-_n(\beta)^d}{D^+_n(\beta)^d}\;.
$$
$$
f^+_n(p,\beta)=\frac{1}{nD^-_n(\beta)^d} \sum_{i=1}^n \frac{1}{\theta_{\bx_i}(p)}K\left(\frac{d_g(p,\bx_i)}{D^+_n(\beta)}\right)=f^+_n(p,D^+_n(\beta)^d)\frac{D^+_n(\beta)^d}{D^-_n(\beta)^d}\;.
$$
where  $D_n^-(\beta)=\beta^{1/2d}h_n$,  $D_n^+(\beta)=\beta^{-1/2d}h_n$  and $h_n^d={k_n}/{(n\lambda({\cal V}_1)f(p))}$ with $\lambda({\cal V}_1)$ denote the Lebesgue measure of the ball in $\real^d$ with radius $r$ centered at the origin.  Note that
\begin{eqnarray}
 \dst\sup_{p\in {{M}_0}}|f^-_n(p,\beta)-\beta f(p)|\convpp 0 \mbox{ and } \dst\sup_{p\in {{M}_0}}|f^+_n(p,\beta)-\beta^{-1} f(p)|\convpp 0. \label{uno}
\end{eqnarray}
 For all $0<\beta<1$ and $\varepsilon>0$ we define
\begin{eqnarray*}
S^-_n(\beta,\varepsilon)&=&\{w:\;\dst\sup_{p\in {{M}_0}}|f^-_n(p,\beta)- f(p)|<{\varepsilon}\;\},\\
S^+_n(\beta,\varepsilon)&=&\{w:\;\dst\sup_{p\in {{M}_0}}|f^+_n(p,\beta)- f(p)|<{\varepsilon}\;\},\\
S_n(\varepsilon)&=&\{w:\;\dst\sup_{p\in {{M}_0}}|\wf_n(p)- f(p)|<{\varepsilon}\;\},\\
A_n(\beta)&=&\{  f^-_n(p,\beta)\leq \wf_n(p)\leq f^+_n(p,\beta)\}
\end{eqnarray*}
Then, $A_n(\beta)\cap S^-_n(\beta,\varepsilon)\cap S^+_n(\beta,\varepsilon)\subset S_n(\varepsilon)$. Let $A=\sup_{p\in M_0}f(p)$. For    $0<\varepsilon <3A/2$  and  $\beta_{\varepsilon}=1-\frac{\varepsilon}{3A}$ consider the following sets
\begin{eqnarray*}
G_n(\varepsilon)&=&\left\{w:\;D_n^-(\beta_{\varepsilon})\leq {\zeta_n(p)} \leq D_n^+(\beta_{\varepsilon})\; \mbox{ for all } p \in M_0 \right\}\\
G^-_n(\varepsilon)&=&\{ \dst\sup_{p\in {{M}_0}}|f^-_n(p,\beta_{\varepsilon})-\beta_{\varepsilon} f(p)|<\frac{\varepsilon}{3}\}\\
G^+_n(\varepsilon)&=&\{ \dst\sup_{p\in {{M}_0}}|f^+_n(p,\beta_{\varepsilon})-\beta_{\varepsilon}^{-1} f(p)|<\frac{\varepsilon}{3}\}.
\end{eqnarray*}
Then we have that $G_n(\varepsilon)\subset A_n(\beta_{\varepsilon})$, $G^-_n(\varepsilon)\subset S^-_n(\beta_{\varepsilon},\varepsilon)$ and $G^+_n(\varepsilon)\subset S^+_n(\beta_{\varepsilon},\varepsilon)$. Therefore,
$G_n(\varepsilon)\cap G^-_n(\varepsilon)\cap G^+_n(\varepsilon)\subset S_n(\varepsilon)$.

On the other hand, using that  $\lim_{r\to 0}{V(B_{r}(p))}/{r^d\mu({\cal V}_1)}=1$, where $V(B_{r}(p))$ denotes the volume of the geodesic ball centered at $p$ with radius $r$ (see Gray and Vanhecke (1979)) and  similar arguments those considered in Devroye and Wagner (1977), we get that
$$
\sup_{p\in M_0}\left|\frac{k_n}{n\lambda({\cal V}_1)f(p)H^d_n(p)}-1\right|\convpp 0.
$$
Recall that $inj_gM>0$ and  $H_n^d(p)\convpp 0$. Then for straightforward calculations we obtained that $ \sup_{p\in M_0}\left|\frac{k_n}{n\lambda({\cal V}_1)f(p)\zeta^d_n(p)}-1\right|\convpp 0.$ Thus, $I_{G^c_n(\varepsilon)}\convpp 0$ and   (\ref{uno}) imply that $I_{S^c_n(\varepsilon)}\convpp 0$. \square

\subsection*{Proof of Theorem \ref{distribucion}.2.}

A Taylor expansion of second order gives
\begin{eqnarray*}
\sqrt{k_n}\left\{\frac{1}{n\zeta_n^d(p)}{\dst\sum_{j=1}^n\dst\frac{1}{\theta_{\bx_j}(p)}K\left(\dst\frac{d_g(p,\bx_j)}{\zeta_n(p)}\right)}-f(p)\right\}=A_n+B_n+C_n
\end{eqnarray*}
where
\begin{eqnarray*}
A_n&=& (h^d_n/\zeta_n^d(p))\sqrt{k_n}\left\{\frac{1}{nh_n^d}{\dst\sum_{j=1}^n\dst\frac{1}{\theta_{\bx_j}(p)}
K\left(\dst\frac{d_g(p,\bx_j)}{h_n}\right)}-f(p)\right\},\\
B_n&=&\sqrt{k_n}((h^d_n/\zeta_n^d(p))-1)\left\{f(p)+\frac{[(h_n/\zeta_n(p))-1]h_n^d}
{[(h^d_n/\zeta_n^d(p))-1]\zeta_n^d(p)}\;\;\dst\frac{1}{nh_n^d}{\dst\sum_{j=1}^n
\dst\frac{1}{\theta_{\bx_j}(p)}K_1\left(\dst\frac{d_g(p,\bx_j)}{\zeta_n(p)}\right)}\right\}  \\
and \\C_n&=& \sqrt{k_n}((h^d_n/\zeta_n^d(p))-1)\frac{[(h_n/\zeta_n(p))-1]^2}{2[(h^d_n/\zeta_n^d(p))-1]}\;\;\dst\frac{1}{n\zeta^d_n(p)}{\dst\sum_{j=1}^n\dst\frac{1}{\theta_{\bx_j}(p)}K_2\left(\dst\frac{d_g(p,\bx_j)}{\xi_n}\right)}[\xi_n/h_n]^2
\end{eqnarray*}
with  $h_n^d={k_n}/{nf(p)\lambda({\cal V}_1)}$ and $\min(h_n,\zeta_n)\leq \xi_n\leq \max(h_n,\zeta_n)$. Note that $H6$ implies that $h_n$ satisfies the necessary hypothesis given in Theorem 4.1 in Rodriguez and Henry (2009), in particular
$$\sqrt{nh_n^{d+4}}\to \beta^{\frac{d+4}d}(f(p)\lambda({\cal V}_1))^{-\frac{d+4}{2d}}.$$
By the Theorem and the fact that $h_n/\zeta_n(p)\convprob 1$, we obtain that $A_n$ converges to a normal distribution  with  mean $b(p)$ and variance $\sigma^2(p)$.  Therefore it is enough to show that $B_n$ and $C_n$ converges to zero in probability.

Note that $\frac{(h_n/H_n(p))-1}{(h^d_n/\zeta_n^d(p))-1}\convprob d^{-1}$ and by similar arguments those considered  in Theorem 3.1 in  Pelletier (2005) and Remark \ref{distribucion}.1. we get that

$$\frac{1}{nh_n^d}{\dst\sum_{j=1}^n\dst\frac{1}{\theta_{\bx_j}(p)}K_1\left(\dst\frac{d_g(p,\bx_j)}{\zeta_n(p)}\right)}\convprob \int K_1(\bu)d\bu f(p)=-d\; f(p).$$  Therefore, by Theorem \ref{distribucion}.3., we obtain that $B_n\convprob 0$. As $\xi_n/h_n$ converges to one in probability, in order to concluded the proof, it remains to prove that $$\dst\frac{1}{n\zeta^d_n(p)} \dst\sum_{j=1}^n\dst\frac{1}{\theta_{\bx_j}(p)}|K_2\left({d_g(p,\bx_j)}/{\xi_n}\right)|$$
is bounded in probability.

By $H7$, there exits $r>0$ such that $|t|^{d+1}|K_2(t)|\leq 1 $ if $|t|\geq r$. Let  $C_r=(-r,r)$, then we have that
\begin{eqnarray*}
\dst\frac{1}{n\zeta^d_n(p)} \dst\sum_{j=1}^n\dst\frac{1}{\theta_{\bx_j}(p)}\left|K_2\left(\dst\frac{d_g(p,\bx_j)}{\xi_n}\right)\right|
&\!\!\leq&\!\!\!\!
\dst\frac{\sup_{|t|\leq r}|K_2(t)|}{n\zeta^d_n(p)} \dst\sum_{j=1}^n\dst\frac{1}{\theta_{\bx_j}(p)}\dst I_{C_r}\!\left(\frac{d_g(p,\bx_j)}{\xi_n}\right)\\
&\!\!+&\!\!\!\!
\dst\frac{1}{n\zeta^d_n(p)} \dst\sum_{j=1}^n\dst\frac{1}{\theta_{\bx_j}(p)}I_{C_r^c}\!\left(\frac{d_g(p,\bx_j)}{\xi_n}\right)\left|\dst\frac{d_g(p,\bx_j)}{\xi_n}\right|^{-(d+1)}\\
\end{eqnarray*}
As $\min(h_n,\zeta_n(p))\leq \xi_n\leq \max(h_n,\zeta_n(p))=\widetilde{\xi}_n$ it follows that
\begin{eqnarray*}
\dst\frac{1}{n\zeta^d_n(p)}\!\!\!\!&&\!\!\!\dst\sum_{j=1}^n\dst\frac{1}{\theta_{\bx_j}(p)}\left|K_2\left(\dst\frac{d_g(p,\bx_j)}{\xi_n}\right)\right|\leq\\
&\leq&
\dst \left(\frac{h_n}{\zeta_n(p)}\right)^d {\sup_{|t|\leq r}|K_2(t)|}\; \frac{1}{nh^d_n}
\dst\sum_{j=1}^n\dst\frac{1}{\theta_{\bx_j}(p)}I_{C_r}\!\left(\frac{d_g(p,\bx_j)}{h_n}\right)\\
&\!\!+&\!\!\!\!
\dst{\sup_{|t|\leq r}|K_2(t)|}\;\frac 1{n\zeta^d_n(p)} \dst\sum_{j=1}^n\dst\frac{1}{\theta_{\bx_j}(p)}I_{C_r}\!\left(\frac{d_g(p,\bx_j)}{\zeta_n(p)}\right)\\
&\!\!+&\!\!\!\!
\dst\left(\frac{h_n}{\zeta_n(p)}\right)^d\frac{1}{nh^d_n} \dst\sum_{j=1}^n\dst\frac{1}{\theta_{\bx_j}(p)}I_{C_r^c}\!\left(\frac{d_g(p,\bx_j)}{h_n}\right)\left|\dst\frac{d_g(p,\bx_j)}{h_n}\right|^{-(d+1)}\left|\dst\frac{\widetilde{\xi}_n}{h_n}\right|^{(d+1)}\\
&\!\!+&\!\!\!\!
\dst\frac{1}{n\zeta^d_n(p)} \dst\sum_{j=1}^n\dst\frac{1}{\theta_{\bx_j}(p)}I_{C_r^c}\!\left(\frac{d_g(p,\bx_j)}{\zeta_n(p)}\right)\left|\dst\frac{d_g(p,\bx_j)}{\zeta_n(p)}\right|^{-(d+1)}\left|\dst\frac{\widetilde{\xi}_n}{\zeta_n(p)}\right|^{(d+1)}\\
&=& C_{n1}+C_{n2}+C_{n3}+C_{n4}.
\end{eqnarray*}
By similar arguments those considered  in Theorem 3.1 in  Pelletier (2005), we have that $C_{n1}\convprob f(p) \int I_{C_r}(s)ds$ and $C_{n3}\convprob f(p)\int I_{C_r^c}(s)|s|^{-(d+1)} ds$.

Finally, let  $A_n^\varepsilon=\{(1-\varepsilon)h_n\leq\zeta_n\leq(1+\varepsilon)h_n\}$ for $0\leq \varepsilon \leq 1.$ Then for $n$ large enough $P(A_n^\varepsilon)>1-\varepsilon$ and in $A_n^\varepsilon$ we have that
\begin{eqnarray*}
I_{C_r}\left(\frac{d_g(\bx_j,p)}{\zeta_n(p)}\right)&\leq& I_{C_r}\left(\frac{d_g(\bx_j,p)}{(1+\varepsilon)h_n}\right),\\
I_{C^c_r}\left(\frac{d_g(\bx_j,p)}{\zeta_n(p)}\right)\left|\frac{d_g(\bx_j,p)}{\zeta_n(p)}\right|^{-(d+1)}&\leq & I_{C^c_r}\left(\frac{d_g(\bx_j,p)}{(1-\varepsilon )h_n}\right)\left|\frac{d_g(\bx_j,p)}{(1-\varepsilon )h_n}\right|^{-(d+1)}\left|\frac{\zeta_n(p)}{(1-\varepsilon )h_n}\right|^{(d+1)}.
\end{eqnarray*}
This fact and analogous arguments those considered in  Theorem 3.1 in  Pelletier (2005), allow to conclude the proof.\square

\subsection*{Proof of Theorem \ref{distribucion}.3.}
Denote $b_n={h_n^d}/(1+zk_n^{-1/2})$, then $$P(\sqrt{k_n}({h_n^d}/{\zeta_n^d}-1)\leq z)=P(\zeta_n^d\geq b_n)=P(H_n^d\geq b_n,\; inj_gM^d\geq b_n).$$ As $b_n\to 0$ and $inj_gM>0$, there exists $n_0$ such that for all $n\geq n_0$ we have that $$P(H_n^d\geq b_n,\; inj_gM^d\geq b_n)=P(H_n^d\geq b_n).$$ Let $Z_i$ such that $Z_i=1$ when $d_{g}(p,\bx_i)\leq b^{1/d}_n$ and $Z_i=0$ elsewhere. Thus,  we have that  $P(H_n^d\geq b_n)=P(\sum_{i=1}^nZ_i\leq k_n)$.
Let  $q_n=P(d_{g}(p,\bx_i)\leq b^{1/d}_n)$. Note that $q_n\to 0$ and $nq_n\to \infty$ as $n\to \infty$, therefore
$$P\left(\sum_{i=1}^nZ_i\leq k_n\right)=P\left(\frac 1{\sqrt{nq_n}}\sum_{i=1}^n(Z_i-E(Z_i))\leq \frac 1{\sqrt{nq_n}}(k_n-nq_n)\right).$$
 Using the Lindeberg Central Limit Theorem we easily obtain that ${(nq_n)}^{-1/2}\sum_{i=1}^n(Z_i-E(Z_i))$ is asymptotically normal with mean zero and variance one. Hence, it is enough to show that ${(nq_n)}^{-1/2}(k_n-nq_n)\convprob z+b_1(p)$.

Denote by $\mu_n=n\dst\int_{B_{b_n^{1/d}}(p)}(f(q)-f(p))d\nu_g(q)$. Note that $\mu_n=n\;q_n-w_n$ with $w_n=
n\;f(p)V(B_{b_n^{1/d}}(p))$. Thus,
 $$\frac 1{\sqrt{nq_n}}(k_n-nq_n)= w_n^{-1/2}(k_n-w_n)\left(\frac{w_n}{w_n+\mu_n}\right)^{1/2}+\frac {\mu_n}{w_n^{1/2}}\left(\frac{w_n}{w_n+\mu_n}\right)^{1/2}.$$

Let $(B_{b_n^{1/d}}(p),\psi)$ be a coordinate normal system. Then, we note that
$$
\frac{1}{\lambda({ \cal V}_{b_n^{1/d}})}\dst\int_{B_{b_n^{1/d}}(p)}f(q)d\nu_g(q)=\frac{1}{\lambda({\cal V}_{b_n^{1/d}})}\dst\int_{{\cal V}_{b_n^{1/d}}}f\circ \psi^{-1}(\bu)\theta_p\circ\psi^{-1}(\bu)d\bu.
$$
The Lebesgue's Differentiation Theorem and the fact that $\dst\frac{V({ B}_{b_n^{1/d}}(p))}{\lambda({\cal V}_{b_n^{1/d}})}\to 1 $ imply that $\dst\frac{\lambda_n}{w_n}\to 0$. On the other hand, from Gray and  Vanhecke (1979), we have that  $${V({ B}_{r}(p))}={r^d\lambda({\cal V}_{1})}(1-\dst\frac{\tau}{6d+12}r^{2}+O(r^{4})). $$
  Hence, we obtain that
\begin{eqnarray*}
w_n^{-1/2}(k_n-w_n)&=& \frac{w_n^{-1/2}\; k_n \;z\;k_n^{-1/2}}{1+zk_n^{-1/2}}+\dst\frac{w_n^{-1/2} \;\tau b_n^{2/d}\;k_n}{(6d+12)(1+zk_n^{-1/2})}+w_n^{-1/2}\;k_n\;O(b_n^{4/d})\\
&=& A_n+B_n+C_n.
\end{eqnarray*}
It's easy to see that $A_n\to z$ and $
w_n^{-1/2} \;b_n^{2/d}\;k_n= \frac{k_nn^{-1/2}b_n^{2/d-1/2}}{(f(p)\lambda({\cal V}_1))^{-2/d}} \left(\frac{b_n\lambda({\cal V}_1)}{V({ B}_{b_n^{1/d}}(p))}\right)^{1/2},
$ since $H6$ we obtain that $B_n\to {\tau\;\beta^{(d+4)/d}}/({6d+12})\;(f(p)\mu({\cal V}_1))^{-2/d}$. A similar argument shows that $C_n\to 0$ and therefore we get that $w_n^{-1/2}(k_n-w_n) \to z+\beta^{\frac {d+4}{d}}\frac{\tau}{6d+12}(f(p)\lambda({\cal V}_1))^{-d/2}$.

In order to concluded the proof we will show that $\dst\frac{\mu_n}{w^{1/2}_n} \to  \frac{\beta^{\frac {d+4}{d}}}{(f(p)\lambda({\cal V}_1))^{(d+2)/d}}\int_{{\cal V}_1}u^2_1\;d\bu\; L_1(p)$. We use a second Taylor expansion that leads to,
\begin{eqnarray*}
\dst\int_{B_{b_n^{1/d}}(p)}(f(q)-f(p))d\nu_g(q)&=&\sum_{i=1}^d \frac{\partial f\circ \psi^{-1}}{\partial u_i}|_{u=0}b_n^{1+1/d}\int_{{\cal V}_1} u_i\;\theta_p\circ\psi^{-1}(b_n^{1/d}\bu)\;d\bu\\
&+&\sum_{i,j=1}^d \frac{\partial^2 f\circ \psi^{-1}}{\partial u_i\partial u_j}|_{u=0}b_n^{1+2/d}\int_{{\cal V}_1} u_iu_j\;\theta_p\circ\psi^{-1}(b_n^{1/d}\bu)\;d\bu\\
&+&O(b_n^{1+3/d}).
\end{eqnarray*}
Using again a  Taylor expansion on $\theta_p\circ \psi^{-1}(\cdot)$  at $0$ we have that
\begin{eqnarray*}
\dst\int_{B_{b_n^{1/d}}(p)}(f(q)-f(p))d\nu_g(q)&=&b_n^{1+2/d}\int_{{\cal V}_1} u_1^2\;d\bu \;L_1(p)+O(b_n^{1+3/d})
\end{eqnarray*}
and by $H6$ the theorem follows.\square

\end{document}